\documentclass[11pt, oneside]{article}   	
\usepackage{geometry}                		
\usepackage{graphicx}				
\usepackage{amssymb}

\usepackage{amsmath,amsfonts,amssymb,amsthm}
\usepackage{cite}
\usepackage{graphics,graphicx}
\usepackage{topcapt}

\usepackage{rotating}
\usepackage{color, colortbl}
\usepackage{multirow}

\definecolor{lblue}{rgb}{0.63,0.79, 0.95}
\definecolor{banana}{rgb}{0.98,0.91,0.71}
\definecolor{lgray}{rgb}{0.86,0.86,0.86}
\definecolor{brightturquoise}{rgb}{0.03, 0.91, 0.87}
\definecolor{darkpastelgreen}{rgb}{0.01, 0.75, 0.24}
\definecolor{junebud}{rgb}{0.74, 0.85, 0.34}
\definecolor{mikadoyellow}{rgb}{1.0, 0.77, 0.05}
\definecolor{mistyrose}{rgb}{1.0, 0.89, 0.88}
\definecolor{aliceblue}{rgb}{0.94, 0.97, 1.0}
\definecolor{antiquewhite}{rgb}{0.98, 0.92, 0.84}
\definecolor{grannysmithapple}{rgb}{0.66, 0.89, 0.63}
\definecolor{lavenderblue}{rgb}{0.8, 0.8, 1.0}
\definecolor{lightsalmon}{rgb}{1.0, 0.63, 0.48}
\definecolor{silver}{rgb}{0.75, 0.75, 0.75}
\definecolor{maize}{rgb}{0.98, 0.93, 0.37}
\definecolor{mintgreen}{rgb}{0.6, 1.0, 0.6}
\definecolor{turquoisegreen}{rgb}{0.63, 0.84, 0.71}
\definecolor{titaniumyellow}{rgb}{0.93, 0.9, 0.0}
\definecolor{tearose(rose)}{rgb}{0.96, 0.76, 0.76}
\definecolor{richbrilliantlavender}{rgb}{0.95, 0.65, 1.0}

\usepackage{epstopdf}

\usepackage{hyperref} 

\hypersetup{colorlinks=true, linkcolor=blue,  anchorcolor=blue,  
citecolor=blue, filecolor=blue, menucolor=blue, pagecolor=blue,  
urlcolor=blue,pdftitle={ACmain}}

\usepackage{cleveref}

\title{Developing Workforce with Mathematical Modeling Skills}
\author{Ariel Cintr\'{o}n-Arias \thanks{Department of Mathematics \& Statistics,
Box 70663,
East Tennessee State University,
Johnson City, TN, 37614-0663, USA (\texttt{cintronarias@etsu.edu}).}
\and Ryan Nivens \thanks{Department of Curriculum \& Instruction,
Box 70684,
East Tennessee State University,
Johnson City, TN,  37614-1709, USA.}
\and Anant Godbole \thanks{Department of Mathematics \& Statistics,
Box 70663,
East Tennessee State University,
Johnson City, TN, 37614-0663, USA.}
\and Calvin B. Purvis \thanks{University Career Services,
Box 70718,
East Tennessee State University,
Johnson City, TN, 37614-0718, USA.}
}

\date{}	

\begin{document}
\maketitle

\begin{abstract}
Mathematicians have traditionally been a select group of academics that produce high-impact 
ideas allowing substantial results in several fields of science. Throughout the past 
35 years, undergraduates enrolling in mathematics or statistics have represented a nearly 
constant rate of approximately 1\% of bachelor degrees awarded in the United States. Even 
within STEM majors, mathematics or statistics only constitute about 6\% of undergraduate degrees 
awarded nationally. However, the need for STEM professionals continues to grow and the 
list of needed occupational skills rests heavily in foundational concepts of 
mathematical modeling curricula, where the interplay of data, computer simulation 
and underlying theoretical frameworks takes center stage.  It is not viable to expect a majority of 
these STEM undergraduates to pursue a double-major that includes mathematics.
Here we present our solution, some early results of implementation, and a vision for possible nationwide adoption.  
\end{abstract}

{\bf Keywords:}  Mathematical Modeling, Mathematics Minor, Modeling, 
Career Development, Career Discovery, Career Services, Employment, STEM Education.

{\bf Funding:} Portions of this work were funded by grant number DUE-1356397, from the National Science Foundation, 
and index number E21454, from the Department of Mathematics \& Statistics of East Tennessee State University. 

{\bf AMS subject classifications:} 97D10, 97M10, 62P99.

\section{Introduction}\label{intro}

\begin{figure}[t]
   \centering
   \includegraphics[keepaspectratio,width=5in,height=3in]{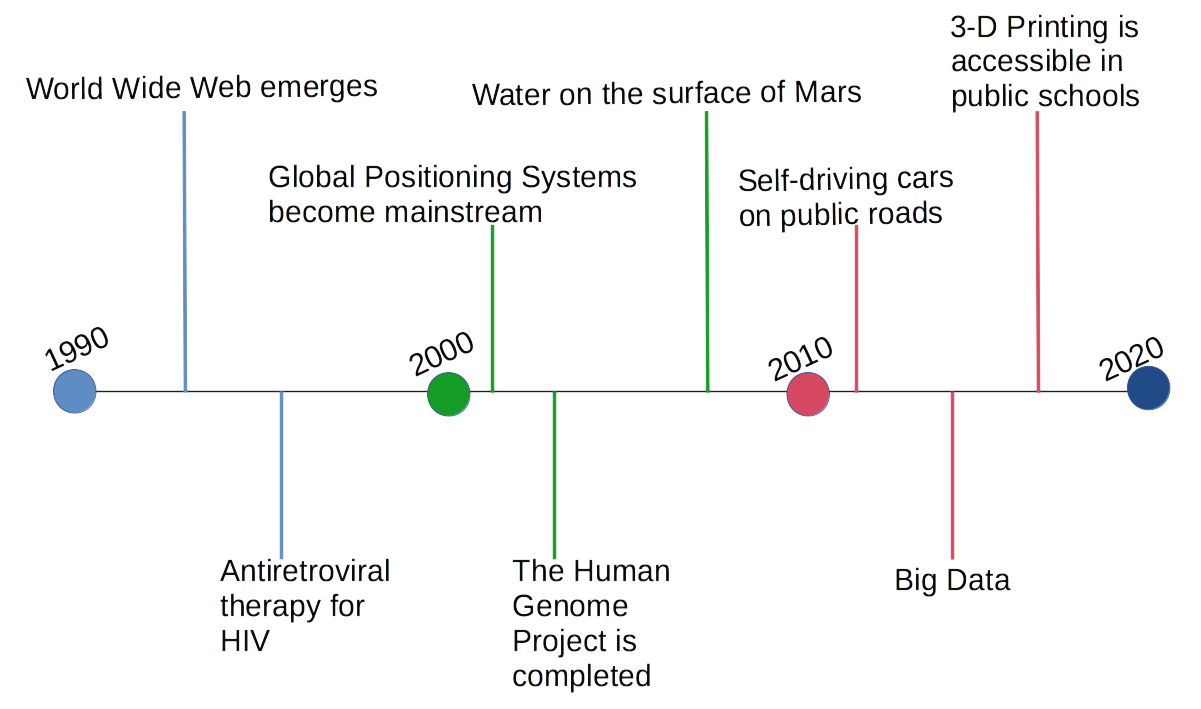} 
   \caption{Timeline of remarkable innovations over three decades.  These accomplishments are interdisciplinary, 
   bringing together professionals in Science, Technology, Engineering and Mathematics (STEM).}
   \label{timelineinnov}
\end{figure}

A period of three decades is short within the timeline of a modern society.  Yet, from 1990 to date we have witnessed scientific and technological innovations 
evolve to become mainstream.  A succinct summary of examples selected by the authors is illustrated in Figure \ref{timelineinnov}.  A common factor 
across these accomplishments is 
interdisciplinary collaboration.  In the latter part of the 20th century, 
one of the worst communicable diseases that circulated
as a global epidemic was caused by the human immunodeficiency 
virus (HIV).  From the mid 1980's through the mid 1990's, interdisciplinary teams of immunologists, virologists, applied mathematicians 
and statisticians eventually succeeded in finding the means to interrupt the life cycle of HIV while extending the life expectancy of infected patients.  An important 
breakthrough came in 1996 from the Theoretical Biology and Biophysics group at Los Alamos National Laboratory, a facility created by 
the United States Department of Energy with an original emphasis on nuclear physics during the era of the Manhattan project, 
when mathematical models describing the interaction between infected cells and virus particles were fitted to viral load 
measurements \cite{perelson96}.   The first estimates from patient data of the viral life cycle duration were incredibly insightful in 
re-designing clinical trials that eventually evolved into successful HIV treatment.

One could argue that two contributing factors of the Big Data revolution in the early 21st century include the establishment of the World Wide Web and 
the completion of the Human Genome Project during the 1990's and the 2000's, respectively.  More and more, professionals in Science, Technology, Engineering, 
and Mathematics (STEM) benefit from extended workforce preparation that includes proficiencies in areas related but outside of their 
primary specialty field.  Examples would include:  
an immunologist who is well-versed in statistical experimental design;
an engineer who is seasoned in data visualization; and
a computer programmer who is competent in molecular biology.  
Integrated STEM is the key to success; the shift from silos in STEM (i.e., the four areas taught and learned separately) to integrated STEM \cite{wang11} is 
the paradigm variation that we believe is fundamental to tackle the next generation of societal problems. 

This article is organized as follows. 
We discuss the nationwide production of STEM professionals in \Cref{stemproduction}, followed by selected STEM occupations in \Cref{stemoccupations}.
\Cref{stemed} conveys how the modeling process relates to the scientific method, STEM education efforts at a local level, 
and modeling across the curriculum at a national scale. Moreover, \Cref{bpmm} summarizes selected educational resources for
mathematical modeling instruction. Our proposed minor in mathematical modeling is discussed in \Cref{modelingminor},
followed by preliminary results, in \Cref{etsusstem},
of a pilot implementation held at East Tennessee State University (home institution of the authors). Concluding remarks are included in \Cref{conclusion}.

\section{STEM Bachelor's Degrees Awarded}\label{stemproduction}

We argue here that undergraduate mathematics programs
across the country are not at capacity in terms of enrollment and graduation rates.
However, substantial growth can be sustained in applied mathematics instruction,
if suitable target populations are brought together with a common foundation, such as mathematical modeling.

\begin{figure}[t]
\begin{center}
\includegraphics[height = 3in,width = 5.5in]{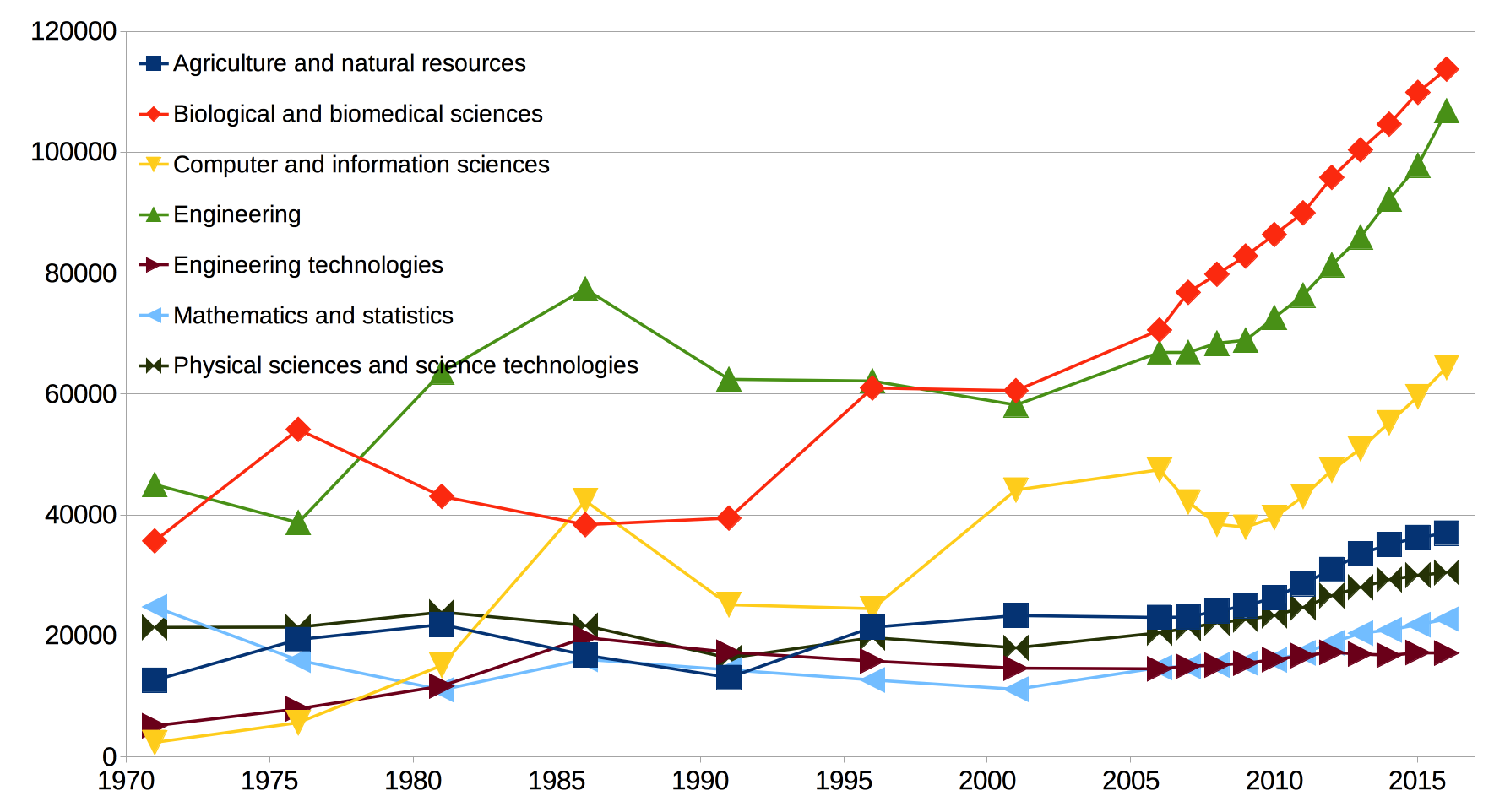}
\includegraphics[height = 3in,width = 5.5in]{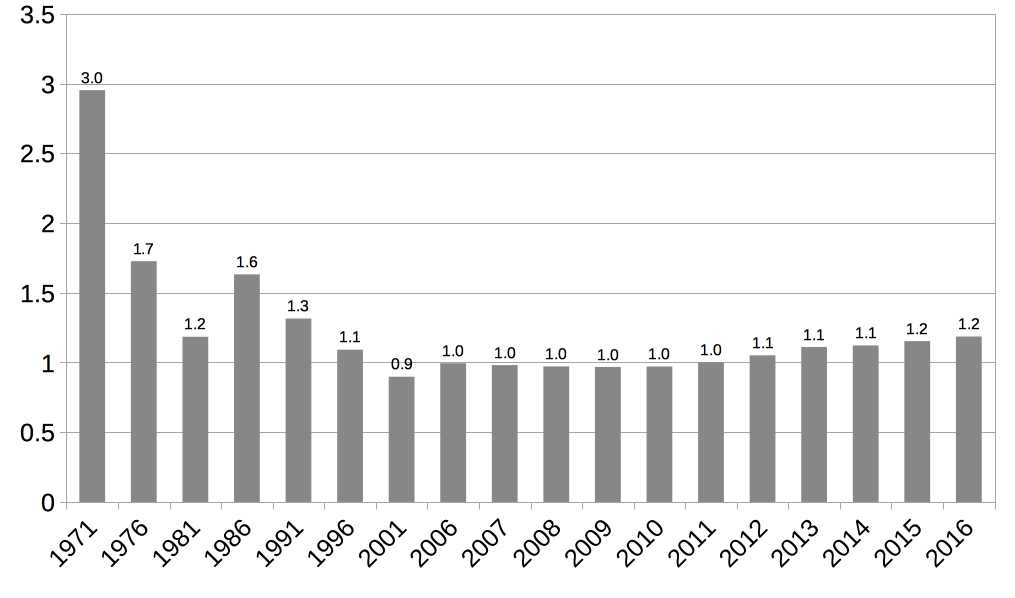}
\caption{The top panel depicts the
number of STEM Bachelor's degrees awarded in the United States, from 1971 through 2016 \cite{nces}.
Seven STEM areas are identified: Agriculture and natural resources (squares); Biological and biomedical sciences (diamonds); 
Computer and information sciences (downward triangles); Engineering (upward triangles); 
Engineering technologies (right triangles); Mathematics and statistics (left triangles); 
Physical sciences and science technologies (double triangles).
The lower panel displays longitudinal percentages of undergraduate mathematics degrees awarded
in the United States, across all disciplines.}
\label{bachelorsstem}
\end{center}
\end{figure}

Longitudinal measurements for the number of Bachelor's degrees in STEM are displayed in Figure \ref{bachelorsstem} (top panel),
with five-year steps from 1971 through 2006 and then in annual steps afterwards.
These counts in Figure \ref{bachelorsstem} were collected by the United States (U.S.) Department of Education \cite{nces} from institutions with
Title IV federal financial aid programs.  

In 1971, the total number of people earning undergraduate degrees, across all disciplines
in the U.S., added up to approximately eight hundred thousand.  By 2016, this total was just under two million \cite{nces}.
On average, 18\% of undergraduate degrees nationwide are awarded in STEM disciplines.  The degrees 
quantified in the top panel of Figure \ref{bachelorsstem} are classified in seven areas (see caption), where 
\emph{Agriculture and natural resources} include agriculture, agriculture operations
\& related sciences, and natural resources \& conservation.  Additionally, \emph{Engineering technologies}
comprise engineering technologies \& engineering-related fields, construction trades, and mechanic \& repair technologies. 
Moreover, \emph{Physical sciences and science technologies} encompass chemistry, physics, astronomy, 
earth sciences, ocean sciences, and atmospheric sciences.

In Figure \ref{timelineinnov}, we considered selected innovations during 1990--2020.  Figure \ref{bachelorsstem} includes data,
in the preceding period 1971--1986, as a baseline.  Growth
is observed in most STEM areas (i.e., positive average rate of change) during 1990--2016.  
Specifically, in 1990 through 2000 there were three (out of seven) STEM areas experiencing growth, 
while throughout 2000--2010 there were five areas, followed by
six areas from 2010 to 2016.  Clearly, biology, engineering, and computing underwent substantial growth
half-way through the second decade of the 21st century.

\clearpage

\clearpage

Nearly 6\% of undergraduate degrees in all STEM areas are in mathematics and/or statistics \cite{nces}, which would naturally supply a pipeline of \emph{``quants"}.
The latter are often referred to as ``the rocket scientists"
of the financial world, they are people prone to mathematical abstraction with substantial understanding of complex mathematical models, who
easily integrate probabilistic and mathematical skills with computing.  

Across all disciplines in the U.S., quants are one-percenters.  This is documented in the lower panel of Figure \ref{bachelorsstem},
where the percentage of individuals earning undergraduate degrees in mathematics and/or statistics is displayed versus time.
For the better part of 35 years, between 1981 and 2016, a steady
state (or horizontal asymptote) is observed, ranging from 0.9\% to 1.6\% with a 5-year average by 2016 equal to 1.1\% \cite{nces}.  According to the 
U.S. Department of Labor \cite{onet}, members of this one percent join the workforce with some of the following job titles:
\begin{itemize}
	\item {Psychometric Consultant}
	\item {Demographer}
	\item {Trend Investigator}
	\item {Agent-Based Modeler}
	\item {Cryptographic Vulnerability Analyst}
	\item {Emerging Solutions Executive}
	\item {Simulation Modeling Engineer}
\end{itemize}

Increasing the number of students majoring in either mathematics or statistics
remains a supreme challenge.  Even though such gains would alter the 1\% across all disciplines nationwide
and the 6\% within STEM, they are expected to be modest or nominal at best, thus implying insufficient
room for growth.  On the other hand, efforts to increase the number of mathematics minors have favorable
odds. Since nearly 99\% of individuals nationwide do not complete a mathematics major \cite{nces}, 
then many of them can be recruited to complete a minor with essential applied mathematics skills.
In the next section,
we review some facts about samples of STEM occupations as reported by the 
U.S. Department of Labor \cite{onet,fastgrowth}.

\section{Selected STEM Occupations}\label{stemoccupations}

The Bureau of Labor Statistics forecasts substantial expected growth in six quantitative STEM occupations \cite{fastgrowth, onet}.
Table \ref{hotskills} reports the projected relative percent change in workforce size,
over the decade 2022 through 2032, for
Computer and Information Research Scientists,
Data Scientists, Information Security Analysts, Operations Research Analysts, Physicists,
and Statisticians.

All the occupations in Table \ref{hotskills} are expected to exceed the average (across all occupations) expected growth of 3\%,
during a decade severely influenced by the challenges of the COVID-19 pandemic. Specific occupation profiles 
can be retrieved in \cite{onet,careeros}, where details such as common duties, typical wages, 
and routine credentials are available.

Table \ref{hotskills} also lists
selected technology skills, which the U.S. Department of Labor
refers to as ``\emph{hot technologies}" \cite{onethot}. Such designation follows from frequently appearing as
requirements in employer job postings. Python and R can be considered the underlying programming
languages driving much of the current innovations in data science and artificial intelligence. 
Structured query language (SQL) is utilized by
several relational database management systems, where
relations among tables storing data are derived from shared fields or attributes (i.e., columns).

The quantitative STEM occupations in Table \ref{hotskills} have undergraduate curricula
that commonly require college mathematics readiness, which serves as a foundation in the
mathematical modeling minor that we propose in \Cref{modelingminor}.

\begin{table}[t]
   \centering
   \topcaption{  
  Selected quantitative STEM occupations with expected employment growth,
   during 2022 -- 2032 \cite{fastgrowth, onet}. Estimates of employment relative percent change (second
   column) were calculated with actual employment in 2022 and projected employment in 2032. 
   Most of the selected occupations are aligned with
   technology skills in high-demand as reported by \cite{onethot}, namely:
   Python, R, and structured query language (SQL). }
\resizebox{\textwidth}{!}{
\begin{tabular}{p{5.5cm}|c|c| c|c|}\hline\hline
Occupation 			& Employment Change & Python		& R	& SQL \\ \hline\hline 
Computer \& Information Research Scientists &23\%& \checkmark&\checkmark&\checkmark  \\ \hline
Data Scientists&35\%& \checkmark&\checkmark&\checkmark  \\ \hline
Information Security Analysts&32\%& \checkmark&&\checkmark  \\ \hline
Operations Research Analysts &23\%& \checkmark &\checkmark&\checkmark  \\ \hline
Physicists &5\% & \checkmark& &\checkmark  \\ \hline
Statisticians&32\%& \checkmark&\checkmark& \checkmark \\ \hline
\end{tabular}}
\label{hotskills}
\end{table}

\section{STEM Education \& the Modeling Process} \label{stemed}

\begin{figure}[h]
   \centering
   \includegraphics[keepaspectratio,width=\textwidth]{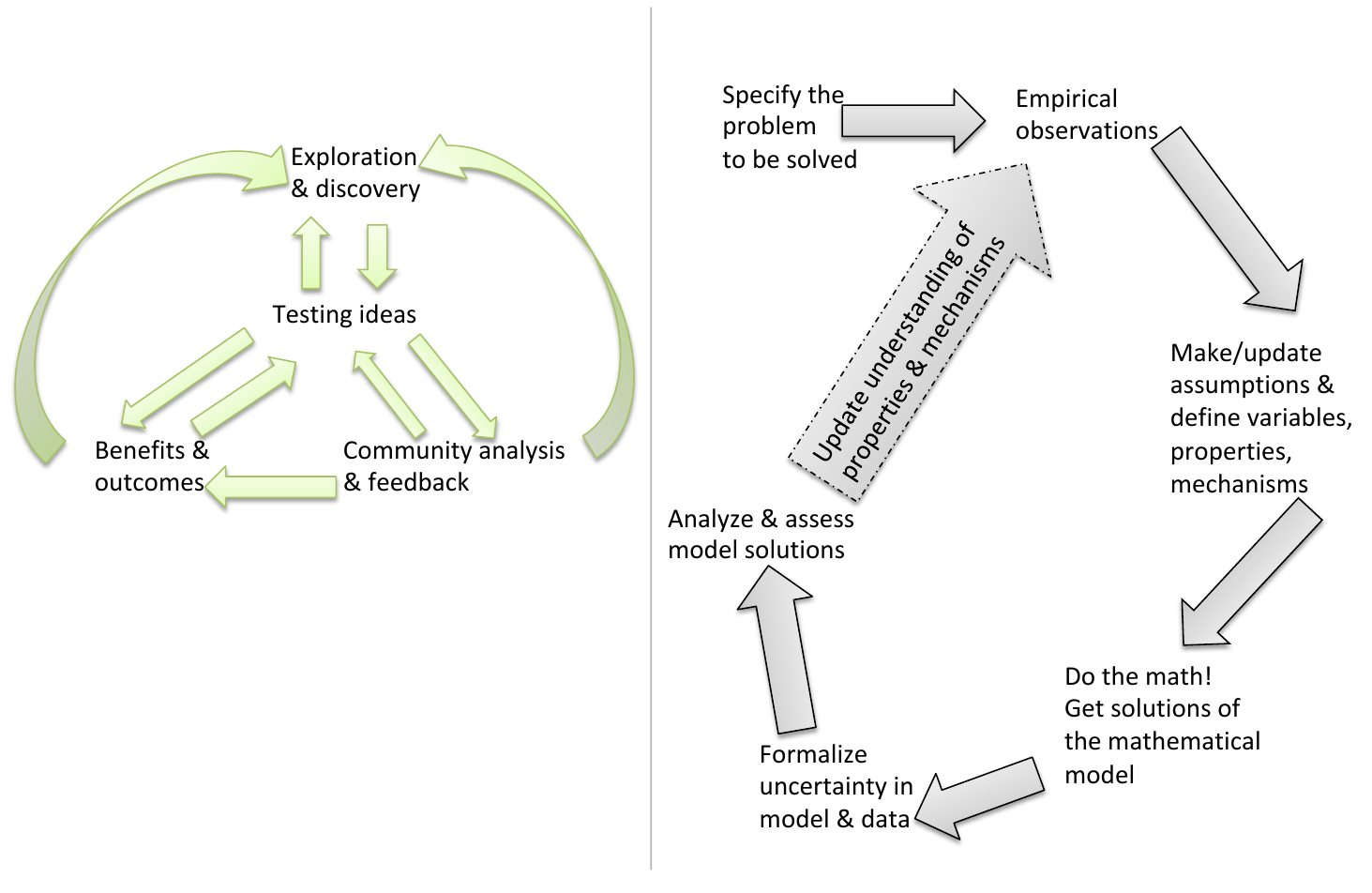} 
   \caption{Schematics for the process of science \cite{scienceprocess}
   and the modeling process \cite{bankstran09,GAIMME}
   are displayed on the left and right panel, respectively.}
   \label{procscimod}
\end{figure}

There are similarities between the process of science and the mathematical modeling process.
The process of science, also commonly referred to as \emph{``the scientific method"}, includes several stages or steps that interactively
evolve into what we know as science.  For many, such a process is a means to do research yielding scientific theories and 
eventually scientific laws.  As explained in \cite{scienceprocess}:
``science is both a body of knowledge and the process for building that knowledge".  The left panel of 
Figure \ref{procscimod} summarizes the process of science.  Testing ideas, including hypotheses, is
central for knowledge generation emerging out of experiments and detailed observations.

The mathematical modeling process is represented by a schematic in the right panel of 
Figure \ref{procscimod}.  Much like the process of science, modeling has various steps that
are expected to be iterated repeatedly,
as described in \cite{bankstran09,GAIMME,mod2012,mod2014,ugappmath14}.
In general, mathematical modeling involves approximations of phenomena in the real-world
that enable analysis and predictions to further understanding  \cite{GAIMME}.  

High-quality STEM education involves an integrated approach \cite{johnson13,nrc02}.  Weaving together science, technology, engineering, and 
mathematics improves student achievement \cite{beckerpark11}, reflects the integrated nature of STEM professions \cite{wang11},
enables deeper understanding of science \cite{nrc02}, and has been observed in successful STEM schools \cite{laforce16}. 
We discuss next a few STEM education efforts implemented in Tennessee, the home state of the authors, during 1990 -- 2020 (the timeline depicted in Figure \ref{timelineinnov}).

Tennessee adopted a value-added assessment system in 1993, long before other states. In March 2010, 
using STEM education as a key component of their grant application, Tennessee was awarded Race to the Top Phase 1 funding
by the U.S. Department of Education.  Then, the Tennessee Department of Education (TDOE), in partnership with
the Battelle Memorial Institute, launched the Tennessee STEM Innovation Network (TSIN), 
committed to developing high-quality STEM curriculum and instruction. Subject to the mantra ``Kindergarten to jobs", TSIN aims to better prepare Tennesseans for the 
STEM workforce by utilizing regional STEM innovation hubs and STEM designated schools located across the state \cite{tnsin}. Under the direction of the STEM 
Leadership Council formed in 2014, and in partnership with the TDOE, the Tennessee STEM Strategic Plan was published in 2016 \cite{tnstem}. 
The plan provides recommendations on full integration of STEM in K-12 education, targeting four priorities to drive the integration of 
mathematics and science standards with broader STEM-related focuses. In 2017, TDOE and TSIN collaboratively designed the \emph{STEM School Designation}. 
Resources were created defining the attributes necessary for a K-12 school to create a comprehensive STEM learning environment. The designation 
has high standards for integration and collaboration, presenting STEM as integrated, rather than insular disciplines, while 
preparing students with needed workforce skills in Tennessee.
An increasing number of K-12 schools are applying for and receiving this designation, 
thus reshaping how STEM education is delivered statewide.  However, across higher education institutions in Tennessee (and in most of the U.S.), 
students currently engage in STEM disciplines separately. Their STEM experience is in silos of science, technology, engineering, or mathematics, 
even though this is very different from how they will experience these fields in the workforce.  We claim that a mathematical modeling minor
may enact the existing national/state recommendations for integrated teaching/learning of interdisciplinary STEM. 

The Society for Industrial and Applied Mathematics (SIAM) made substantial advances during the decade 2010 -- 2020 while
improving STEM education, by fostering mathematical modeling across the K-16 curriculum. Strategy, methodology, and materials were
conveyed in pivotal reports published by SIAM in 2012 \cite{mod2012}, 2014 \cite{mod2014}, and 2019 (second edition) \cite{GAIMME}.

The first report, \emph{Modeling Across the Curriculum} \cite{mod2012}, documents recommendations by subject matter experts
during a workshop (funded by the National Science Foundation) with four main pillars: 
(1) Expand modeling in K-12 STEM curriculum;
(2) Develop high school modeling courses;
(3) Develop modeling-based undergraduate curricula, with initial efforts geared towards STEM majors during their first year;
(4) Develop a repository of instructional materials, such as: lesson plans, articles, books, websites, videos, and problems with solutions.

The second report, \emph{Modeling Across the Curriculum II} \cite{mod2014}, originated from a follow-up workshop 
with the following aims: (1) Incorporate mathematical modeling for early grades; 
(2) Develop high school curricular materials; 
(3) Create two modeling pathways for undergraduate STEM education.

The GAIMME (Guidelines for Assessment and Instruction of Mathematical Modeling Education) report has
educators as the intended primary audience \cite{GAIMME}. It stipulates:
``Unlike a course on traditional topics in which students may listen to lectures 
or follow examples to solve a set of problems, students in a 
mathematical modeling course tackle big, messy, open-ended problems without any textbook examples." \cite{GAIMME}.  In fact, this principle
could be a resonating factor that shifts unengaged students (receiving mathematics instruction that lacks
meaningful applications) into individuals who actively engage with components of the modeling process, because they
understand how it relates to their career goals (and in some cases, to their personal lives) \cite{GAIMME}.
Insights into optimizing the benefits of the GAIMME report
in undergraduate curricula, while further developing modeling skills, can be found in \cite{bliss19}.

\section{Educational Resources in Mathematical Modeling}\label{bpmm}
A primer into the process of mathematical modeling, geared towards students and teachers, is available in
\emph{Math Modeling: Getting Started and Getting Solutions} \cite{bliss14}. The companion handbook to the latter is
\emph{Math Modeling: Computing and Communicating} \cite{bliss18}, which expands from the modeling process into technical 
computing while using software platforms and coding. Moreover, guidelines for higher education instructors are available in
Chapter 4 and Appendices A \& D of the GAIMME report \cite{GAIMME}.

{\bf \underline{Software}.} Traditionally, several tasks involved in the modeling process (right panel of Figure \ref{procscimod})
are implemented with MATLAB, Maple or Mathematica, to mention a few. Here, based on the data reported by the 
U.S. Department of Labor (see Table \ref{hotskills}), and the teaching
history of one of the authors, we suggest to adopt open source software, such as:
Python, R, and/or SQL. An instructional best-practice while adopting Python and R is to deploy Jupyter notebooks with two
main approaches: (1) local install of Anaconda \cite{conda}; (2) web-based interface with either Google Colab \cite{colab}
or CoCalc \cite{cocalc}.
Professional development, free of cost for educators, is available from 
Data Carpentry \cite{datacpy,datacr}, Data Camp \cite{datacampu}, and Data Quest \cite{dataquest}.

{\bf \underline{Datasets}.} We now mention some highlights of three public data repositories that can
be used with instruction of mathematical modeling. First, the
\emph{Systemic Initiative for Modeling Investigations \& Opportunities with Differential Equations (SIMIODE)} \cite{simiode},
is a community of practice where datasets are available together with instructional guides, worksheets, and peer-reviewed content.
Next, the U.S. General Services Administration maintains a public website in compliance with the 
\emph{Open, Public, Electronic and Necessary (OPEN) Government Data Act}, which releases data, tools, and resources
to design data visualizations \cite{datagov}. Lastly, the \emph{UCI Machine Learning Repository} is a collection of more 
than six hundred databases, domain theories, and data generators that are used by 
for the empirical analysis of predictive modeling algorithms \cite{uciml}.

\section{Mathematical Modeling Minor}\label{modelingminor}

\begin{figure}[t]
   \centering
   \includegraphics[keepaspectratio,width=\textwidth]{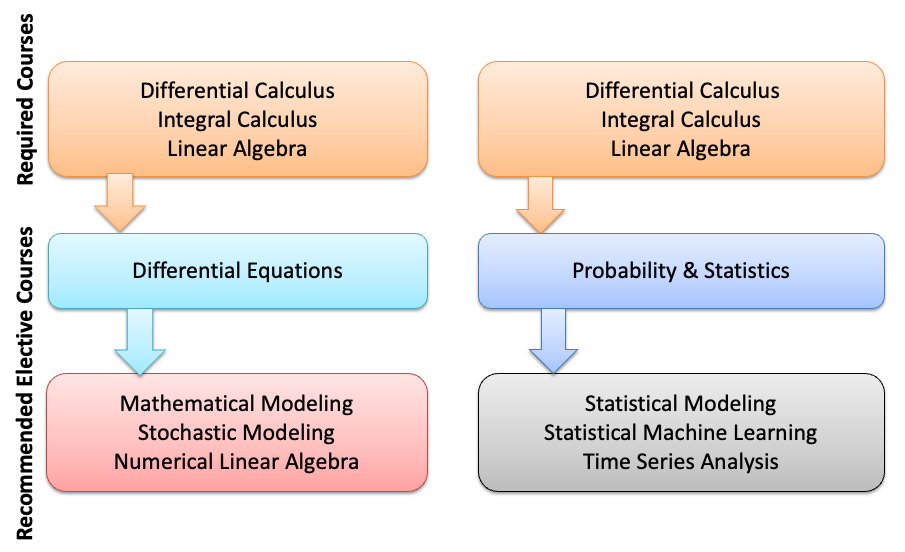}
   \caption{Mathematics minor course sequences resembling two concentrations offered at East Tennessee State University (ETSU). The left-side panel 
	resembles the concentration in Computational Applied Mathematics, while the right-side panel mimics the Statistics concentration.}
   \label{seqminor}
\end{figure}

The home institution of the authors is
East Tennessee State University (ETSU), which offers a Bachelor of Science degree in mathematics 
with four possible concentrations.  Two of them are directly related to mathematical modeling, specifically:  \emph{Computational Applied Mathematics}
and \emph{Statistics}.  Figure \ref{seqminor} illustrates sequences of classes that echo both of these concentrations 
in the context of a mathematics minor, also offered at ETSU. Even though the formal name is generic (i.e., \emph{mathematics minor}) 
the rules in the ETSU catalog are flexible enough to foster plenty of modeling skills while fulfilling the requirements of a minor.

In 2014, the \emph{SIAM Education Committee Report on Undergraduate Degree Programs in Applied Mathematics} \cite{ugappmath14},
provided an outline for an applied mathematics minor. We propose here an updated version with emphasis on modeling, influenced by 
the modeling process (Figure \ref{procscimod}) and technology skills in high-demand (Table \ref{hotskills}):
\begin{itemize}
	\item[] \underline{Required Courses} 
		\begin{itemize}
			\item Calculus sequence
			\item Linear Algebra
			\item Probability and Statistics
		\end{itemize}
	\item[] \underline{Recommended Elective Courses} 
		\begin{itemize}
			\item Coding for Data Manipulation
			\item Differential Equations with Modeling
			\item Numerical Methods for Deep Learning
			\item Computational Linear Algebra
			\item Mathematical Modeling
			\item Statistical Modeling
			\item Stochastic Modeling
			\item Predictive Modeling
			\item Operations Research
			\item Time Series Analysis
		\end{itemize}
\end{itemize}

Every institution of higher education would role out a version of this proposed minor in mathematical modeling
that best fits the expertise of departmental faculty and availability in its course catalog. We recommend
six-to-seven courses, with a total of 18-to-23 credits. Next, we include some remarks about selected courses 
from the list above, mostly informed by our own instructional experiences.

\clearpage

{\bf \underline{Calculus sequence}.} The bare minimum requirement would consist of only one course,
namely: Differential Calculus. In some instances, Integral Calculus and/or Multivariate Calculus 
may already be requirements of a STEM major and yet only a subset of a three-course sequence would suffice 
as minor requirements, allowing more room for recommended elective courses.
Often, calculus readiness is attained
by completing a general education course in either Precalculus or College Algebra. It is important not to overload
the minor with hidden pre-requisites. 

{\bf \underline{Coding for Data Manipulation}.} This course is intended as an introduction to computer programming with 
Python, R, and SQL, within the context of tabular data. Key topics would include: control flow, object-oriented programming, 
grammar of graphics, tidy data, combining tabular data, and reshaping tabular data. Guidelines can be found in
an article by H. Wickham \cite{wickhamt} and textbooks by W. McKinney \cite{mcKinney} and H. Wickham \cite{wickhamg}.

{\bf \underline{Differential Equations with Modeling}.} A two-course calculus sequence is a common pre-requisite.
Materials provided by SIMIODE \cite{simiode} and the textbook by K. Bryan \cite{bryank} are optimal blueprints for this course. 

{\bf \underline{Numerical Methods for Deep Learning}.} All minor required courses may serve as pre-requisites.
We recommend an evolved version of a scientific computing course that 
aligns with applications of machine learning. Key topics may include: Root-finding methods, interpolation and splines,
numerical solvers of differential equations, implementations of neural networks,
Monte Carlo sampling, Latin hypercube methods,
ensemble methods, gradient-based algorithms, and direct search methods. A textbook by B. Despr\'{e}s \cite{depresb}
can provide useful guidelines. 

\begin{table}[t]
   \centering
   \topcaption{
   Two versions of a generic minor in mathematical modeling, based on recommendations by SIAM Education Committee \cite{ugappmath14} and
   skills in high-demand as reported by the U.S. Department of Labor (see Table \ref{hotskills} and references therein). Each version is a sequence of
   seven courses.
   }
\begin{tabular}{|c|c|}\hline\hline
		Minor without Differential Equations& Minor with Differential Equations\\ \hline\hline
		 Differential Calculus &  Differential Calculus \\ \hline 
		 Linear Algebra & Integral Calculus \\ \hline 
		 Probability \& Statistics & Linear Algebra \\ \hline 
		 Coding for Data Manipulation & Probability \& Statistics \\ \hline 
		 Mathematical Modeling & Differential Equations with Modeling \\ \hline 
		 Stochastic Modeling & Time Series Analysis \\ \hline 
		 Computational Linear Algebra& Mathematical Modeling\\ \hline 
\end{tabular}
   \label{2ver}
\end{table}

{\bf \underline{Computational Linear Algebra}.} This course would require Linear Algebra, Probability and Statistics, and
Differential Calculus as pre-requisites. The following topics may be covered: tensor products, matrix factorization, matrix decompositions,
topic modeling with singular value decompositions, principal component analysis, and compressed sensing with robust regression.
Open lessons and Jupyter notebooks can be dowloaded from GitHub \cite{fastainla}. Additionally, further details about 
the underlying Python libraries for this course are available in a textbook by J. Howard and S. Gugger \cite{howard2020},
together with code templates in GitHub \cite{fastai}.

{\bf \underline{Predictive Modeling}.} Pre-requisites may include: Probability and Statistics,
Differential Calculus, and Linear Algebra. Essential topics comprise:
machine learning life cycle, data preprocessing, feature engineering, exploratory data analysis,
overfitting and underfitting, tree-based models, regression models, ensembling models, hyperparameter tuning, and bias mitigation.
Blueprint teaching materials can be obtained with zero cost for students and educators from the
AWS Academy \cite{awsacademy}, IBM SkillsBuild \cite{ibmsb}, and the Cisco Networking Academy \cite{cskills4all}.

Mathematical modeling involves coding skills and data manipulation. Based on our own experience, we know
there are two possible approaches:
(1) requiring a stand alone course (similar to the one described above) where basic skills in control flow, syntax,
and data frame transformations are established; (2) integrating coding skills in some of the required courses and in all of the
recommended elective courses. One of the authors successfully implemented the latter,
by deploying self-paced supplemental instruction with 
interactive smart materials provided at low-cost for students by vendors such as Data Camp \cite{datacampu} and Data Quest \cite{dataquest}.

Another consideration is requiring both Differential Calculus and Integral Calculus, instead of just the former.
In Table \ref{2ver} we display two scenarios, each of which with a total of seven required courses
for a mathematical modeling minor. On the left-side, only Differential Calculus is part of the required courses but
Coding for Data Manipulation is included as a recommended elective course. In the contrast, the right-side
illustrates a scenario where there are four required courses (including a two-course sequence of Calculus), and
three recommended elective courses, encompassing Differential Equations with Modeling.
Also, the right-side of Table \ref{2ver} displays a scenario where the coding skills would be integrated in several courses, 
instead of having a stand alone course.

Undergraduate research is an `` ... investigation conducted by an undergraduate student that makes an original intellectual or creative 
contribution to the discipline" \cite{ugappmath14,cur}. Students completing a mathematical modeling minor will benefit from 
applying their modeling skills within undergraduate research projects that originate from their major discipline.
There are several U.S. federal government agencies that provide funding for undergraduate research, here we mention three of them:
the U.S Department of Energy \cite{suli}, the National Security Agency \cite{nsareu}, and the National Science Foundation \cite{nsfreu}.

\section{Case Study: Preparation of Data-Driven Mathematical Scientists for the Workforce}\label{etsusstem}
Impactful funding for undergraduate mathematical modeling immersion was provided to ETSU by 
a grant from the National Science Foundation, under a program in the  
Division of Undergraduate Education called \emph{Scholarships in 
Science, Technology, Engineering, and Mathematics} (S-STEM).

The formal name of this S-STEM program at ETSU was 
\emph{Preparation of the Data-Driven Mathematical Scientists for the Workforce}, and it ran from
August 1, 2014 through July 31, 2020. In its initial years, this program provided specialty training 
in the form of special seminars. In latter years, it evolved by integrating modeling modules within existing curriculum
housed by the Department of Mathematics \& Statistics at ETSU.

Even though the first goal was to fund as many as 30 students, 
the gross total of funded scholars equals 58, where 43 of them did successfully complete the program.
Longitudinal enrollment in the S-STEM program
had several academic terms when the participants divided almost evenly by gender.
However, from total head counts there were nearly one third female versus two thirds
male scholars.  Most of the participants were recruited from the mathematics major (which has almost an even split by gender)
and the computing major (where females constitute less than one third) at ETSU.

In Table \ref{sstemskills} we summarize some of the hard skills that S-STEM scholars acquired during the program.  
Computing tools, such as, Linux, Python, and R, were first introduced in special topics seminars and would then be adopted into some of the ETSU curriculum.  
Students were introduced to Micro-Array Data Analysis and Predictive Modeling in workshops and tutorials.
All other topics in Table \ref{sstemskills} were delivered as upper-level courses that had integrated tasks in R or Python (promoting
two of the hot technologies in Table \ref{hotskills}).
Modeling was a unifying factor, as noted in the title of the courses Mathematical Modeling, Statistical Modeling, and Stochastic Modeling.  
Due to the complex task of lining up schedules from different majors, it was not possible for all scholars to take exactly the same sequence of courses.  
Instead, they would combine coursework from the sequences depicted in Figure \ref{seqminor}.


\begin{table}[h]
   \centering
   \topcaption{Hard skills developed by S-STEM participants at ETSU.} 
   \resizebox{\textwidth}{!}{	
\begin{tabular}{|m{3.25cm}|m{3.25cm}|m{3.25cm}|m{3.25cm}|} \hline
\cellcolor{white}Linux, Python, R & \cellcolor{white}   Predictive Modeling & \cellcolor{white}Micro-Array Data Analysis& \cellcolor{white}Mathematical Modeling\\ \hline
 \cellcolor{white} Statistical Machine Learning & \cellcolor{white} Applications of Statistics & \cellcolor{white}Numerical Analysis & \cellcolor{white} Statistical Modeling \\  \hline
  \cellcolor{white} Stochastic Modeling & \cellcolor{white}  Probability \& Statistics & \cellcolor{white} Numerical Linear Algebra & \cellcolor{white} Bayesian Analysis \\ \hline
\end{tabular}}
   \label{sstemskills}
\end{table}

Several of the S-STEM scholars were highly motivated to apply for summer programs.  
Part of the motivation was drawn by direct conversations with guest speakers that visited East Tennessee State University, sponsored by the S-STEM program.  
Institutions where scholars obtained paid summer positions (either as an internship or as a research experience for undergraduates) include
universities (e.g., Arizona State University), businesses (e.g., Eastman Chemical Company), and federal government agencies (e.g. National Security Agency).
We argue that participation in summer programs, outside their home institution for most, was a contributing factor to improve the graduation rate of S-STEM scholars.
The four-year graduation rate of S-STEM scholars ranged between 67\% and 91\%, reflecting the positive impact provided by the 
workforce preparation scope of this program. In contrast, ETSU's four-year graduation rate did range from 23\% to 32\%\cite{ipeds}. 

Several S-STEM participants continued graduate studies either in Master's degree 
programs (4 students) or doctoral programs (6 students). S-STEM scholars that did not continue graduate studies 
upon graduation were successful landing employment, 
in positions that require skills with an interplay between data, computer simulations, and mathematics and/or statistics.  
The following are some of their employers in the private sector and federal government:
\begin{itemize}
	\item Clayton Homes
	\item Eastman Chemical Company
	\item National Security Agency
	\item Riparian LLC
	\item U.S. Army Materiel Systems Analysis Activity
\end{itemize}

Recruiting full-time, first-time students into the S-STEM program was a challenge. We kept a public webpage and circulated public service announcements in
a local National Public Radio station,
yet the number of submitted applications was nominal, at best. After the second year, we expanded our strategy by targeting students enrolled in the Calculus sequence and
Linear Algebra, which required developing a relationship with the instructors alongside with an onset of traits that were used to filter potential applicants, by the middle of a semester.
After pools of potential applicants were identified, we placed individual phone calls delivering ``an elevator pitch" of the S-STEM program.

Just taking one course over a summer term can make a substantial difference in two or three subsequent semesters. This was another area of growth for the S-STEM program at ETSU. 
In Summer 2018, the courses Integral Calculus, Differential Equations, and Linear Algebra were offered with online asynchronous delivery, over a period of ten weeks.
S-STEM scholars did attend synchronous supplemental instruction sessions (while using video conferencing). Prior to this summer, 
these courses were only offered in-person, over an intense period of five weeks. Our instructional design targeted adult learners working full-time during the 
summer at remote locations. Thus, completing these courses with remote delivery
in Summer 2018, allowed S-STEM scholars to enroll in upper-level classes (e.g., see Figure \ref{seqminor}) during Fall 2018 and Spring 2019.

Having grant funding (award number DUE-1356397) for STEM majors was essential while boosting
enrollment in a modeling-enhanced mathematics minor at ETSU.
It successfully compounded institutional benefits in five areas: recruitment, retention, internships, graduation,
and employment.  Data-driven mathematics were promoted with modeling at the helm, and they
inspired the ideas that we share here.  In the next section, we close with some concluding remarks.

\section{Final Remarks}\label{conclusion}
Careers requiring skills in data-driven decision making domains are
expected to widen, over the current decade. In this article, we discussed six
quantitative STEM careers with favorable employment forecasts
during 2022 -- 2032 (see Table \ref{hotskills}).

Technology skills in high-demand within digital economies were cross-examined against
a sample of quantitative STEM occupations in Table \ref{hotskills}. Such skills were also integrated in the design
of a proposed mathematical modeling minor (see \cref{modelingminor}), partly inspired by workforce upskilling programs driven by digital badges and micro-credentials.

Mathematical modeling aligns well with principles of computational thinking, defined by
the International Society for Technology in Education (ISTE) \cite{istect}
as ``a problem-solving process that includes:
logically organizing and analyzing data; representing data through abstractions such as
models and simulations; automating solutions through algorithmic thinking".
One of the authors is an ISTE Certified Educator with a successful history of aligning computational thinking
standards with college mathematics curriculum.

As of this writing, in a COVID era with digital economies of scale, 
we believe a mathematical modeling minor is helpful in workforce development
for digital infrastructures centered around data-driven technologies.



\end{document}